\documentclass[12pt]{article}%22.12.2016
\usepackage{amsmath}
\usepackage{amssymb}

\usepackage{dsfont}
\usepackage{cite}
\newsymbol \blackbox 1004
\newcommand{\eh}{\hfill}\newlength{\sperr}

\def\nn{\nonumber}

\def\la{\lambda}
\def\om{\omega}
\def\Om{\Omega}
\def\t{\theta}

\def\wh{\widehat}
\def\wt{\widetilde}
\def\ov{\overline}

\def\prt{\partial}

\def\BC{{\mathbb C}}
\def\BR{{\mathbb R}}

\def\cla{{\mathcal A}}

\def\clu{{\mathcal U}}

\newcommand{\E}{\mathrm{e}}
\newcommand{\I}{\mathrm{i}}

\newtheorem{Pa}{Paper}[section]
\newtheorem{Tm}[Pa]{{\bf Theorem}}

\newtheorem{Rk}[Pa]{{\bf Remark}}

\newtheorem{Pn}[Pa]{{\bf Proposition}}

\title{Inversion of the convolution operators \\ on a rectangular}

\author{Alexander Sakhnovich}

\date{}
\begin{document}
\maketitle

\begin{abstract}    Using simultaneously two operator identities, we consider the inversion of the convolution operators
on a rectangular. The structure of the inverse operators and of some corresponding forms, which are important in
signal processing, is derived.
\end{abstract}



\section{Introduction}\label{Intro}
\setcounter{equation}{0}

We consider a convolution operator $S$ (or, equivalently, an operator $S$ with difference kernel) of the form
 \begin{align}\label{F1}&
Sf=\frac{\prt}{\prt x_1}\frac{\prt}{\prt x_2}\int_{\Om}s(x-t)f(t)dt \quad \big(x=(x_1, \, x_2), \,\, t=(t_1, \, t_2)\big), 
 \end{align}
where
\begin{align}\label{F2} &
 \Om=\{x: \, 0<x_1<\om_1, \, 0<x_2<\om_2\}, \quad s(x-t)=s(x_1-t_1, \, x_2-t_2). 
  \end{align}
  We assume that $s(x) \in L^2(\wt \Om)$, where 
  $\wt \Om=\{x: \, |x_1|<\om_1, \, |x_2|<\om_2\},$
  and so the integrals in \eqref{F1} are well-defined for $f \in L^2(\Om)$. Moreover, we assume that the right-hand side of \eqref{F1}
  is well-defined and that $S$ is bounded in $L^2(\Om)$.
 
 In one-dimensional case, the inversion of convolution operators is connected with the names of N. Wiener, E. Hopf, M.G. Krein,
 I.C. Gohberg, V.A. Ambartsumian, V. V. Sobolev, L.A. Sakhnovich and many other mathematicians and applied scientists.
 The inversion of convolution operators on a semi-axis (of Wiener-Hopf operators) was studied in various papers including the brilliant works
 \cite{GoKr, Kre} (see also, e.g.,  \cite{CrL, doT, Go, KreS} and references therein). The situation with the inversion of convolution operators $S$
 on a finite interval is more complicated and essentially different from the case of semi-axis. One of the first works on the subject was again
 written by M.G. Krein \cite{Kre0}. Then, a procedure to recover the operator $S^{-1}$, from its action on two functions only, was published
 \cite{SaL73}. Further developments as well as various applications and  references can be found in \cite{SaL80, SaL15}.
 The method of operator identities, which was introduced in \cite{SaL73} (see also \cite{SaL68}), may be successfully used for the
 inversion of various other structured operators. See, for instance, \cite{FKRS, SaA73} and a more detailed discussion with references
 in \cite[Appendix D]{SaSaR}.
%%%%%%%%%%%%%%%%%

According to \cite{SaL15}, ``an important part in the theory of the equations with a difference kernel (on the interval)
is played by the equations with a special right-hand side":
$$Sf(x, \la)=\E^{\I  \la x}.$$
The reflection coefficient  in the light scattering theory coincides with the function
 \begin{align}\label{F3}&
\rho(\la,\mu)=\big(S^{-1}\E^{\I  \la x}, \E^{\I  {\ov \mu}x}\big)_{\om} \quad (\la, \mu \, \in \BC), 
 \end{align} 
where $(\cdot, \cdot)_{\om}$ stands for the scalar product in $L^2(0,\om)$ and $\BC$ denotes the complex plane.
The structure of $\rho$ derived in \cite{SaL80} (see also \cite{SaL15}) generalizes a well-known astrophysics
 result from \cite{Amb, Iv, Sob}.
 
The inversion of the convolution operators and the structure of the inverse operators in the two-dimensional case is much more complicated
than in the  one-dimensional case (similar to some other inversion and interpolation problems). In spite of several books (see, e.g., \cite{Nap, Vol}) and many
papers on multidimensional convolution operators, the problem to take simultaneously into account the dependence of the integral kernel of the operator
on the difference of the arguments (with respect to both variables) remained unsolved. This problem is essential in applications (e.g., in astrophysics
and signal processing \cite{JeVa, NiB}) and in factorization and interpolation theories. 

In our paper we express the $\rho$-function
 \begin{align}\label{F4}&
\rho(\la,\mu)=\int_0^{\om_1}\int_0^{\om_2} \E^{-\I  { \mu}x}\big(S^{-1}\E^{\I  \la x}\big)(x_1,x_2)dx_1dx_2  \,\, (\mu x=\mu_1x_1+\mu_2x_2),
\end{align} 
 where $x\in \Om$ and $\la,\mu\in \BC^2$, in terms of a bounded operator acting from
 $L^2_2(0,\om_i)$ into  $L^2_2(0,\om_k)$ $(i\not=k)$. Such an expression is a direct analog
 of \cite[(1.3.12)]{SaL15}, where $\rho$ in one-dimensional case (see \eqref{F3}) is expressed via four functions depending on one variable.
 Using the structure of $\rho$, we derive the structure of the operator, which is inverse to a convolution operator on $\Om$.

Here, the notation ${\mathbf 1}$ stands for the function which identically equals $1$ on the corresponding set, and col means column
(e.g., ${\mathrm{col}}\begin{bmatrix}a & b\end{bmatrix}=\begin{bmatrix}a \\ b\end{bmatrix}$).
As usual, $\BR$ denotes the real axis and $\BC$ denotes the complex plane. 
%The notation $\cld$ stands for the differentiationin the planar Lebesgue measure.
 The class of bounded operators acting from the Hilbert space ${\bf H}_1$ into the Hilbert space ${\bf H}_2$ is denoted by $B({\bf H}_1,{\bf H}_2)$.
 When ${\bf H}_1={\bf H}_2={\bf H}$, we write $B({\bf H})$ instead of $B({\bf H},{\bf H})$.

 %%%%%%%%%%%%%%%%%%%%%%%%%%%%%%%%%%%

%%%%%%%%%%%%%%%%%%%%%%%%%%%%%%%%%%%%%%%%%%%%%%%%%%%%%

\section{Structures of the $\rho$-function \\ and of the inverse operator}\label{rho}
\setcounter{equation}{0}
\subsection{Bounded operators on a rectangular}
Similarly to the proof of \cite[Theorem 1.1.1]{SaL15} on the representation of the operators from
$B\big(L^2(0,\ell)\big)$  (one-dimensional case), one can show that operators 
$Q\in B\big(L^2(\Om)\big)$ admit representation
\begin{align}\label{Fd1} &
Qf= \frac{\prt}{\prt x_1}\frac{\prt}{\prt x_2}\int_{\Om}q(x,t)f(t)dt,
  \end{align}
where the integrals $\int_{\Om}q(x,t)f(t)dt$ are generating functions of  absolutely continuous (with respect to the planar Lebesgue measure)
charges on $\Om$, and for each fixed $x\in L^2(\Om)$ (and varying $t$)  we have $q(x,t)\in L^2(\Om)$. More precisely, we have 
\begin{align}\label{Fd1'} &
q(x,t)=\ov{\big(Q^*\chi_x\big)(t)}, 
  \end{align}
where
\begin{align} \nn & \chi_x(t)=1 \,\, {\mathrm{for}} \,\,  t\in \clu_x=\{t:\, 0<t_1< x_1, \, 0<t_2<x_2\}, \\
\nn &
 \chi_x(t)=0 \,\, {\mathrm{for}} \,\, t\in \big(\Om\backslash \clu_x\big).
\end{align}
We consider a slightly more general than \eqref{F1} class of bounded convolution operators, namely,
operators of the form
 \begin{align}\nn
Sf=\frac{\prt}{\prt x_1}\frac{\prt}{\prt x_2}\int_{\Om}q(x,t)f(t)dt, \quad q(x,t)=& s(x-t)-s(x-t_1,-t_2)
\\ & \label{Fd2}
-s(-t_1,x_2-t_2)+s(-t),
 \end{align}
where the  expression $-s(x-t_1,-t_2)
-s(-t_1,x_2-t_2)+s(-t)$ may be used for normalization purposes and disappears (after the differentiation of the integral in \eqref{Fd2})
in the
cases of comparatively smooth functions $s(x-t)$.

We always assume that $s(x)\in L^2(\wt \Om)$, for $\wt \Om=\{x: \, |x_1|<\om_1, \, |x_2|<\om_2\}$. Without loss of generality we assume also that
\begin{align} & \label{Fd3}
\int_0^{\om_1}s(-t_1,x_2)dt_1=\int_0^{\om_2}s(x_1, -t_2)dt_2=0.
\end{align}
Note that if \eqref{Fd3} does not hold, we may substitute $s(x)$ with 
$$\wh s(x)=s(x)-\frac{1}{\om_1}\int_0^{\om_1}s(-t_1,x_2)dt_1  -\frac{1}{\om_2}\int_0^{\om_2}s(x_1,-t_2)dt_2+
\frac{1}{\om_1\om_2}\int_{\Om}s(-t)dt,$$  
and the function $\wh s$ satisfies \eqref{Fd3}.
%%%%%%%%%%%%%%%%%%%%%%%%%%%%%%%%%%%%%%%%%%%%%%%%%%%%%%%%
%%%%%%%%%%%%%%%%%%%%%%%%%%%%%%%%%%%%%%%%%%%%%%%%%%%%%%%%%
\subsection{Structure of the $\rho$-function}
The $\rho$-function given by \eqref{F4} admits representation
 \begin{align}\nn &
\rho(\la, \mu)=\Big((I-\mu_1A_1^*)^{-1}(I-\mu_2A_2^*)^{-1}S^{-1}(I-\la_2A_2)^{-1}(I-\la_1A_1)^{-1}{\mathbf 1},  {\mathbf 1}\Big)_{\Om}, 
 \end{align}
 where $(\cdot , \cdot)_{\Om}$ is the scalar product in $L^2(\Om)$, ${\mathbf 1}$ is the function which identically equals $1$ (equals $1$ on $\Om$
 in the formula above), and
  the operators $A_k\in B\big(L^2(\Om)\big)\,$  are given by  
\begin{align}\label{F6}&
A_1f=\I \int_0^{x_1}f(t_1, x_2)dt_1, \quad A_2f=\I \int_0^{x_2}f(x_1, t_2)dt_2.
 \end{align}
Our approach is based on the simultaneous usage of two operator identities:
\begin{align}\label{F5}&
A_kS-SA_k^*=\I(M_{1k}M_{2k}+M_{3k}M_{4k}) \qquad (k=1,2).
 \end{align} 
More precisely, we use the next easy proposition. (Recall that we always assume that \eqref{Fd3} is valid and that $s(x) \in L^2(\wt \Om)$.)
\begin{Pn}\label{OpId} Let $S\in B\big(L^2(\Om)\big)$ be an operator
of the form \eqref{Fd2}.

Then, the operator identities 
\eqref{F5}, where $M_{1k}, \, M_{3k}\in B\big(L^2(0, \om_i), \, L^2(\Om)\big)$, $\,\,M_{2k}, \, M_{4k}\in B\big(L^2(\Om), \, L^2(0, \om_i)\big)$
$(i\not= k)$,
\begin{align}\label{F7}&
\big(M_{11}f\big)(x)=\frac{\prt}{\prt x_2}\int_0^{\om_2}s(x_1,x_2-t_2)f(t_2)dt_2, \\
\label{F8}&
\big(M_{12}f\big)(x)=\frac{\prt}{\prt x_1}\int_0^{\om_1}s(x_1-t_1, x_2)f(t_1)dt_1;
\\ \label{F9}&
\big(M_{21}f\big)(x_2)=\int_0^{\om_1}f(t_1,x_2)dt_1, \quad \big(M_{22}f\big)(x_1)=\int_0^{\om_2}f(x_1,t_2)dt_2;
\\ \label{F10}&
\big(M_{31}f\big)(x)=f(x_2), \quad \big(M_{41}f\big)(x_2)=-\frac{\prt}{\prt x_2}\int_{\Om}s(-t_1,x_2-t_2)f(t)dt;
\\ \label{F11}&
\big(M_{32}f\big)(x)=f(x_1), \quad \big(M_{42}f\big)(x_1)=-\frac{\prt}{\prt x_1}\int_{\Om}s(x_1-t_1, -t_2)f(t)dt;
 \end{align}
are valid.
\end{Pn}
The operator identities \eqref{F5} may be rewritten in the more traditional form:
\begin{align}\label{F5'}&
A_kS-SA_k^*=\I\Pi_k \wh \Pi_k, \quad \Pi_k:= \begin{bmatrix} M_{1k} & M_{3k}\end{bmatrix},
\quad \wh \Pi_k:= \begin{bmatrix} M_{2k} \\ M_{4k}\end{bmatrix},
 \end{align} 
where $k=1,2$. Clearly, the operators $A_k^*$ in \eqref{F5'} have the form
\begin{align}\label{F12}&
A_1^*f=-\I \int_{x_1}^{\om_1}f(t_1, x_2)dt_1, \quad A_2^*f=\I \int_{x_2}^{\om_2}f(x_1, t_2)dt_2.
 \end{align}

We note that the integral kernels of the operators $M_{1k}$ and $M_{4k}$ given in \eqref{F7}, \eqref{F8} and in \eqref{F10}, \eqref{F11}, respectively,
depend on the difference of a one of two variables. Introducing integration operators $\cla_k\in B\big(L^2(0, \om_k)\big)$:
\begin{align}\label{F13}&
\cla_1f=\I \int_0^{x_1}f(t_1)dt_1, \quad \cla_2f=\I \int_0^{x_2}f(t_2)dt_2,
 \end{align}
and operators  $K_{1i}\in  B\big(L^2(0, \om_k), \, L^2(0,\om_i)\big)$ $(k\not=i)$, $K_{2i}\in  B\big(\BC, \, L^2(0,\om_i)\big)$, $K_4\in  B\big(L^2(\Om), \BC)\big)$:
\begin{align}\label{F14}&
K_{11}f=- \int_0^{\om_2}s(x_1,-t_2)f(t_2)dt_2, \quad K_{2i}1={\mathbf 1}, 
\\ \label{F15}&
K_{12}f=- \int_0^{\om_1}s(-t_1,x_2)f(t_1)dt_1, \quad K_4=\int_{\Om}s(-t)f(t)dt,
 \end{align}
one easily obtains, for instance, operator identities for $M_{4k}$ (see the proposition below).
\begin{Pn}\label{OpIdfM} Let the conditions of Proposition \ref{OpId} hold. Then the following operator identities
are valid:
\begin{align}\label{F16}&
\cla_i M_{4k}-M_{4k}A_i^*=\I(K_{1i}M_{2i}+K_{2i}K_{4}) \qquad (i,k=1,2, \quad i\not=k).
 \end{align} 
\end{Pn}
Assuming that $S$ has a bounded inverse operator, we express $\rho(\la,\mu)$ in terms of the operators $g_{ik}$ $(i\not=k)$:
\begin{align}\label{F17}&
g_{ik}=\begin{bmatrix} K_{3i} \\ K_{1i}\end{bmatrix}\begin{bmatrix} I & 0\end{bmatrix}- \wh \Pi_k S^{-1} \Pi_i ,
\quad g_{ik}\in B\big(L^2_2(0,\om_k), \, L^2_2(0,\om_i)\big),
 \end{align} 
where the operators $ \Pi_i$ and  $\wh \Pi_k$ are given in \eqref{F5'}, the operators 
$K_{1i}$ are defined in \eqref{F14} and \eqref{F15}, and
\begin{align}\label{F18}&
K_{3i}=\int_0^{\om_k}\cdot \, dt_k\in B\big(L^2(0, \om_k), \, L^2(0, \om_i)\big), \quad i \not=k.
 \end{align} 
The operators $g_{ik}$ and $g_{ki}$ are connected via the relation
\begin{align}\label{F19}&
g_{ki}=-U_k J_k g_{ik}^*J_i U_i; \quad \big(U_i f\big)(x_i)= \ov{f(\om_i -x_i)}, \quad J_i=\I \begin{bmatrix}0 & -I \\ I &0  \end{bmatrix},
 \end{align} 
where $J_i$ is acting in $L^2_2(0, \om_i)\big)$ and $U_i$ is acting (depending on the context) in $L^2_2(0,\om_i)$ or in
$L^2(0,\om_i)$.
\begin{Tm} \label{Tmrho}Let $S\in B\big(L^2(\Om)\big)$ be a convolution operator
of the form \eqref{Fd2} $($where $s(x) \in L^2(\wt \Om))$.  Assume that $S$ has a bounded inverse operator.

Then the $\rho$-function introduced in \eqref{F4} has the forms
\begin{align}\label{F20}&
\rho(\la, \mu)=\frac{1}{\mu_k-\la_k}\int_0^{\om_i}\E^{-\I\om \mu}\big(J_i U_i \psi_i(\mu,x_i)\big)^*\psi_i(\la, x_i)d x_i ,
 \end{align} 
where $i=1$ or $i=2$ $(k=1,2; \,\, k\not=i)$, $\psi_i(\la,x_i)\in L^2_2(0,\om_i)$,
\begin{align}\label{F21}&
\psi(\la,x)=\begin{bmatrix}\psi_1(\la,x_1)\\ \psi_2(\la,x_2)\end{bmatrix}=\theta(\la)G(\la)^{-1}{\mathrm{col}}\begin{bmatrix}0 & {\mathbf 1} & 0 & {\mathbf 1}\end{bmatrix},
\\ \label{F22}&
\sbox0{$\begin{matrix}I-\la_1 \cla_1 & 0\\ 0 &I-\la_1\cla_1\end{matrix}$}
\sbox1{$\begin{matrix}I-\la_2 \cla_2 & 0\\ 0 &I-\la_2\cla_2\end{matrix}$}
G(\la):=\left[
\begin{array}{c|c}
\usebox{0}&\makebox[\wd0]{\large $\I \la_2 g_{12}$}\\
\hline
  \vphantom{\usebox{0}}\makebox[\wd0]{\large $\I \la_1 g_{21}$}&\usebox{1}
\end{array}
\right], 
 \end{align} 
{\rm{col}} in \eqref{F21} means column, and  the scalar function $\t(\la)$ is the unique function of the form
\begin{align}\label{F23}&
\t(\la)=1+\la_1\la_2\int_{\Om}\E^{\I\la(\om-x)}h(x)d x, \quad h\in L^2(\Om),
 \end{align}
 such that the expression on the right-hand side of \eqref{F21} does not have poles or zeros. 
 
 Here the operators $\cla_k$ are introduced in \eqref{F13}, the operators $U_i$ and $J_i$ are introduced in \eqref{F19},
 and the operators $g_{ik}$ are expressed  $($in terms of $S$ and $S^{-1})$ in \eqref{F17}. Moreover $g_{12}$ is expressed via $g_{21}$ and vice
 versa in \eqref{F19}, and so the representations \eqref{F20} are  completely determined by $g_{12}$ or, equivalently, by $g_{21}$.
\end{Tm}
We have also a precise formula for $h(x)$ in \eqref{F23}:
\begin{align}\label{F24}&
h=S^{-1}y, \quad y(x):=s(x_1-\om_1,x_2-\om_2).
 \end{align}
\begin{Rk} Clearly, $\rho(\la, \mu)$ determines $($see \eqref{F25}$)$ the inverse operator $S^{-1}$. According to Theorem \ref{Tmrho}, $\rho(\la, \mu)$ $($and so the operator $S^{-1})$
is determined by $g_{ik}$. 

Recall that the operator $S$ of the form \eqref{Fd2} is determined by its integral kernel  or, equivalently, by the four
functions $s(x), \, s(-x), \, s(x_1,-x_2)$ and $s(-x_1, x_2)$, where $x\in \Om$. Since $g_{ik}\in B\big(L^2_2(0,\om_k), \, L^2_2(0,\om_i)\big)$, the integral
kernel of the operator  $g_{ik}$
is also determined by some four functions on $\Om$.

The considerations above confirm a {\rm heuristic principle} which, for the case of convolution operators on the interval, was formulated
and proved in \cite{SaL80} $($see also \cite{SaL73}$)$. This principle states that the amount of information $($i.e., the number of functions$)$ which
determines $S$ coincides with the minimal amount of information which is necessary to exactly construct $S^{-1}$.

\end{Rk}

\subsection{Structure of the inverse operator}
In view of \eqref{F4}, \eqref{Fd1} and \eqref{Fd1'}, the operator $T=S^{-1}$ admits representation
\begin{align}\label{F25}
\big(Tf\big)(x)=&-\frac{1}{16 \pi^4}\frac{\prt}{\prt x_1}\frac{\prt}{\prt x_2}  \int_{\Om}\left({\mathrm{l.i.m.}}_{r\to \infty}\int_{\clu_{r}}\E^{-\I \la t}\right.
\\ \nn 
& \times  \left.  \int_{\BR^2}(\mu_1\mu_2)^{-1} \rho(\la,\mu)
(\E^{\I \mu_1 x_1}-1)(\E^{\I \mu_2 x_2}-1)d\mu d \la 
\right) f(t)d t,
 \end{align}
where l.i.m. stands for the limit in $L^2$-norm and 
$$\clu_r=\{\la :\, |\la_i|<r \quad (i=1,2)\}.$$
We note that the boundedness of the convolution operator $S$ and of $S^{-1}$ yields the boundedness of the operator $\Gamma$
determined by the relation 
\begin{align}\label{F26}&
\Gamma  \E^{\I \la x}=-\I \begin{bmatrix} \la_2 I & 0 \\ 0 & \la_1 I\end{bmatrix}^{-1}\left(\psi(\la) - 
{\mathrm{col}}\begin{bmatrix} 0 & \E^{\I \la_1 x_1} & 0 & \E^{\I \la_2 x_2} \end{bmatrix} \right),
\end{align}
and acting from $L^2(\Om)$ into $L^2_2(0,\om_1) \oplus L^2_2(0,\om_2)$.

Now, we formulate an inverse result.
\begin{Tm} \label{TmStrT} Let a given operator $g_{ik}$ $($where $i=1,\, k=2$ or $i=2,\, k=1)$ belong to
$B\big(L^2_2(0,\om_k), L^2_2(0,\om_i)\big)$. Assume that the operator $T$, which is determined by $g_{ik}$
via formula \eqref{F25} and via the procedure to construct $\rho(\la,\mu)$ from Theorem \ref{Tmrho}, is bounded.
Assume that the operator $\Gamma$  $($which is determined by $g_{ik}$
via equality \eqref{F26}, relation\eqref{F19} and formulas from Theorem \ref{Tmrho}, namely, formulas
\eqref{F21}-\eqref{F23}$)$  
 is bounded as well.
 
 Then, if the inverse operator $S=T^{-1}$ exists and is bounded, this operator $S$ is a convolution
 operator, that is, $S$ admits representation \eqref{Fd2}.
\end{Tm}

\bigskip

\noindent{\bf Acknowledgments.} The author is grateful to B. Kirstein and L. Sakhnovich for useful discussions.
 {The research   was supported by the
Austrian Science Fund (FWF) under Grant  No. P29177.}
%%%%%%%%%%%%%%%%%%%%%%%%%%%%%%%%%%%%%%%%%%%%%%
%%%%%%%%%%%%%%%%%%%%%%%%%%%%%%%%%%%%%%%%%%%%%%%

%%%%%%%%%%%%%%%%%%%%%%%%

\end{document}